\documentclass[a4paper,11pt]{amsart}
\usepackage{amssymb,amsmath}
\usepackage[xdvi]{graphicx}
\usepackage{sect}
\usepackage{hyperref}

\def\le{\leqslant}

\def\C{\mathbb{C}}
\def\R{\mathbb{R}}
\def\conj{conj}
\newcommand{\sminus}{\smallsetminus}

\swapnumbers
\newtheorem{theom}{\itshape Theorem}[section]

\newtheorem{cor1}[theom]{\itshape Corollary 1}
\newtheorem{cor2}[theom]{\itshape Corollary 2}
\newtheorem{lem}[theom]{\itshape Lemma}

\theoremstyle{remark}
\newtheorem*{rem}{Remark}

\begin{document}

\title[Encomplexing the writhe]{Encomplexing the writhe}

\author{Oleg Viro}
\dedicatory{Uppsala University, Uppsala, Sweden;\break
POMI, St.~Petersburg, Russia}
\address{\newline
Matematiska institutionen, Box 480, Uppsala
universitet, 751 06, Uppsala, Sweden}

\email{oleg@math.uu.se}



\begin{abstract} For a nonsingular real algebraic curve in
$3$-dimensional projective space or $3$-sphere, a new integer-valued
characteristic is introduced. It is invariant under rigid isotopy
and multiplied by $-1$ under mirror reflections. In a sense, it is
a Vassiliev invariant of degree $1$ and a counterpart of a link
diagram writhe.
\end{abstract}

\maketitle

\section{Introduction}\label{sI}

This paper is a detailed version of my preprint \cite{Viro0},
which was written about five years ago.
Here I do not discuss
results that have appeared since then. I plan to survey them soon
in another paper. The subject is now evolving into a
\textit{real algebraic knot theory}.

This paper is dedicated to the memory of my teacher Vladimir
Abramovich Rokhlin. It was V.~A.~Rokhlin, who suggested to me,
a long time ago, in 1977, to develop a theory of real algebraic
knots. He suggested this as a topic for my second dissertation
(after PhD, like habilitation). Following this suggestion, I moved
then from knot theory and low-dimensional topology to the
topology of real algebraic varieties. However, in the topology of
real algebraic varieties, problems on spatial surfaces and plane
curves were more pressing than problems on spatial curves, and my
second dissertation defended in 1983 was devoted to the constructions
of real algebraic plane curves and spatial surfaces with
prescribed topology.

The change in the topic occured mainly because I managed to
obtain decent results in another direction, on plane curves. There
was also a less respectable reason: I failed to relate the
traditional techniques of classical knot theory to real
algebraic knots. One of the obstacles was a phenomenon which
became the initial point of this paper. A large part of the
traditional techniques in knot theory uses plane knot
diagrams, i.e., projections of knots to the plane. The projection of an
algebraic curve is algebraic, and one could try to apply results
on plane real algebraic curves. However, the projection contains
extra real points, which do not correspond to real points of the
knot. These points are discussed below. In the seventies they
ruined my weak attempts to study real algebraic knots. Now they
allow us to detect crucial differences between topological and
real algebraic knots.

I am grateful to Alan Durfee, Tobias Ekholm, and V.~M.~Kharlamov
for stimulating conversations.

The lengthy informal introduction, which follows, is intended to explain
the matter prior to going into details. I cannot resist the temptation
to write in the style of popular mathematics and apologize to the
reader whom this style may irritate.

\subsection{Knot theory and algebraic geometry}\label{sI0}
In classical knot theory, by a link one means a smooth closed
$1$-dimensional submanifold of the $3$-dimensional sphere $S^3$, i.e.,
the union of several disjoint circles smoothly embedded into $S^3$.
A link may be equipped with various additional structures such as
orientation or framing and considered up to various equivalence
relations like smooth (or ambient) isotopy, PL-isotopy, cobordism
or homotopy. See, e.g., \cite{Rolfsen} or \cite{Burde-Zieschang}.

In algebraic geometry classical links naturally appear as links of
singular points of complex plane algebraic curves. Given a singular
point $p$ of a complex plane algebraic curve $C$, the intersection of
$C$ with the boundary of a sufficiently small ball centered at $p$ is
called the \textit{link of the singularity}. It provides a base for a
fruitful interaction between topology and algebraic geometry with
a long history and lots of important results.

Another obvious opportunity for interaction between algebraic
geometry and knot theory is based on the fact that a classical
link may emerge as the set of real points of a real algebraic
curve. This opportunity was completely ignored, besides that a
number of times it was proved that any classical link is
approximated by (and hence isotopic to) the set of real points of
a real algebraic curve. There are two natural directions in which
algebraic geometry and knot theory may interact in the study
of real algebraic links: first, the study of relationships between
invariants which are provided by link theory and algebraic
geometry, second, developing a theory parallel to the classical
link theory, but taking into account the algebraic nature of the
objects. From the viewpoint of this second direction, it is more
natural to consider real algebraic links up to isotopy consisting
of real algebraic links, which belong to the same continuous
family of algebraic curves, rather than up to smooth isotopy in
the class of classical links. I call an isotopy of the former kind
a \textit{rigid isotopy}, following the terminology established by
Rokhlin \cite{R} in a similar study of real algebraic plane
projective curves and the likes (see, e.g., the survey \cite{Viro
New pr.}). Of course, there is a forgetting functor: any real
algebraic link can be regarded as a classical link and a rigid
isotopy as a smooth isotopy. It is interesting to see how much is lost
under that transition.

In this paper I point out a real algebraic link invariant
which is lost. It is unexpectedly simple. In an obvious sense it is a
nontrivial Vassiliev invariant of degree $1$ on the class of real
algebraic knots (recall that a knot is a link consisting of one
component). In classical knot theory the lowest degree of a
nontrivial Vassiliev knot invariant is $2$. Thus there is an essential
difference between classical knot theory and the theory of real
algebraic knots. Of course this difference has a simple topological
explanation: a real algebraic link is more complicated topologically,
besides its set of real points contains the set of complex points
invariant under the complex conjugation and a rigid isotopy induces an
equivariant smooth isotopy of this set.

The invariant of real algebraic links which is defined below is
very similar to the self-linking number of a framed knot. In
\cite{Viro0} I call it also the \textit{self-linking number}. Its
definition looks like a replacement of an elementary definition of
the writhe of a knot diagram, but taking into consideration the
imaginary part of the knot.

\subsection{The word `encomplex'}\label{sI1.5}
Here I propose to change this name (i.e., self-linking number) to
\textit{encomplexed writhe}, and, in general,
since many other characteristics can also be enhanced in a similar
way, I suggest a new verb \textit{encomplex} for similar enhancements by
taking into consideration additional imaginary ingredients.
This would agree with the general usage of the prefix \textit{en-}
which is described in the Oxford Dictionary of Current English as
follows: ``en- \textit{prefix} $\dots$ forming verbs $\dots$ \textbf{1} from
nouns, meaning `put into or on'
(\textit{engulf\textup; entrust\textup; embed}),
\textbf{2} from nouns or adjectives, meaning `bring into the condition
of' (\textit{enslave}) $\dots$''.

The word \textit{complexification} does not seem to be appropriate for
what we do here with the writhe. A complexification of the writhe
should be a complex counterpart for the writhe, it should be a
characteristic of complex objects, while our enhancement of the writhe
is defined only for real objects possessing
complexification.

\subsection{Self-linking and writhe of nonalgebraic knots}\label{sI2}
The linking number is a well-known numerical characteristic of a
pair of disjoint oriented circles embedded in three-dimensional
Euclidean space. Roughly speaking, it measures how many times one
of the circles runs around the other. It is one of the most
classical topological invariants, introduced in the nineteenth
century by Gauss \cite{Gauss}.

In the classical theory, a self-linking number of a knot is
defined if the knot is equipped with an additional structure like
a framing or just a vector field nowhere tangent to the
knot.\footnote{A framing is a pair of
normal vector fields on the knot orthogonal to each other.
There is an obvious construction
that makes a framing from a nontangent vector field and
establishes a one-to-one correspondence between homotopy classes of
framings and nontangent vector fields. The vector fields are more
flexible and relevant to the case.} The self-linking number is the
linking number of the knot oriented somehow and its copy obtained
by a small shift in the direction specified by the vector field.
It does not depend on the orientation, since reversing the
orientation of the knot is compensated by reversing the induced
orientation of its shifted copy. Of course, the self-linking
number depends on the homotopy class of the vector field.

A knot has no natural preferable homotopy class of framings, which
would allow us to speak about a self-linking number of the knot
without a special care on the choice of the
framing.\footnote{Moreover, the self-linking number is used to
define a natural class of framings: namely, the framings with
self-linking number zero.} Some framings appear naturally in
geometric situations. For example, if one fixes a generic
projection of a knot to a plane, the vector field of directions of
the projection appears. The corresponding self-linking number is
called the \textit{writhe} of the knot. However, it depends on the
choice of the projection and changes under isotopy.

The linking number is a Vassiliev invariant of order $1$ of
two-component oriented links. This means that it changes by a
constant (in fact, by $2$) when the link experiences a homotopy with
the generic appearance of an intersection point of the components.
Whether the linking number increases or decreases depends only on
the local picture of orientations near the double point: when it
passes from 
$\vcenter{\hbox{\includegraphics[bb=0 0 16 16,scale=.7,clip]{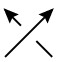}}}$ 
 through
$\vcenter{\hbox{\includegraphics[bb= 0 0 16 16,scale=.7,clip]{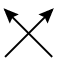}}}$ 
to
$\vcenter{\hbox{\includegraphics[bb= 0 0 16 16,scale=.7,clip]{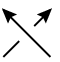}}}$, 
the linking number
increases by $2$. Generalities on Vassiliev invariants see, e.g., in
\cite{V}.

In a sense the linking number is the only Vassiliev invariant of
degree $1$ of two-component oriented links: any Vassiliev invariant
of degree $1$ of two-component oriented links is a linear function
of the linking number. Similarly, the self-linking number is a
Vassiliev invariant of degree $1$ of framed knots (it changes by
$2$ when the knot experiences a homotopy with a generic appearance of
a self-intersection point) and it is the only Vassiliev of degree
$1$ of framed knots in the same sense. The necessity of a framing for
the definition of self-linking number can now be formulated more
rigorously: only constants are Vassiliev invariants of degree $1$ of
(non-framed) knots.

The diagrammatical definition of the writhe, which is imitated
below, runs as follows: for each crossing point of the knot
projection, one defines a \textit{local writhe} equal to $+1$ if
near the point the knot diagram looks like
$\vcenter{\hbox{\includegraphics[bb= 0 0 16 16,scale=.7,clip]{fo3.eps}}}$ 
and $-1$ 
if it looks like
$\vcenter{\hbox{\includegraphics[bb=0 0 16 16,scale=.7,clip]{fo1.eps}}}$. 
Then one sums the local
writhes over all double points of the projection. The sum is the
writhe.

A continuous change of the projection may cause the vanishing of a crossing
point. This happens under the first Reidemeister move shown in the left
hand half of Figure~\ref{f1}. This move changes the writhe by $\pm 1$.

\subsection{Algebraicity enhances the writhe}\label{SI3}
\begin{figure}[ht]
\centerline{\includegraphics[bb= 0 0 352 68,scale=.9,clip]{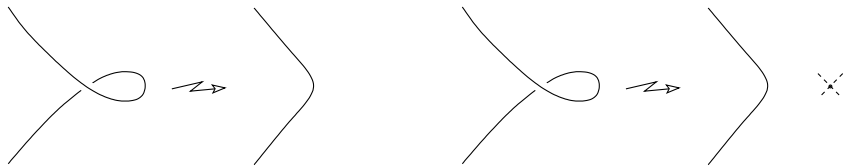}}
\caption{Topological (left) and real algebraic (right) versions of
the first Reidemeister move. At the solitary crossing point, which
is on the right hand side of the picture, the conjugate imaginary
branches are indicated by dashed segments, according to an
outdated tradition of Analytic Geometry.}\label{f1}
\end{figure}

If a link is algebraic, then its projection to a plane is algebraic,
too. A generic projection has only ordinary double points and the
total number of its complex double points is
constant.\footnote{Here by a generic projection we mean a
projection from a generic point. When one says that a generic
projection has some properties, this means that for an open
everywhere dense set of points the projection from any point of
this set has these properties. The whole set of
undesirable points is closed nowhere dense although it depends on
the properties under consideration. A proof is an
easy exercise either on Sard's Lemma, or Bertini's Theorem.} The
number of real double points can vary, but only by an even number.
A real double point cannot turn alone into an imaginary one, as it
seems to happen under the first Reidemeister move. Under an
algebraic version of the first Reidemeister move, the double point
stays in the real domain, but becomes solitary, like the only real
point of the curve $x^2+y^2=0$. The algebraic version of the first
Reidemeister move is shown in the right hand half of Figure~\ref{f1}.

It is not difficult to prove that
the family of spatial curves that realizes this move can be transformed
by a local diffeomorphism to the family of affine curves defined by
the following system of equations
$$
\left\{\begin{aligned}
xz+y&=0,\\
x+z^2+\tau&=0,
\end{aligned}\right.
$$
where $\tau$ is the parameter of the deformation. These are rational
curves, admitting a rational parametrization
$$
\left\{\begin{aligned}
x&=-t^2-\tau,\\
y&=-t(t^2+\tau),\\
z&=-t.
\end{aligned}\right.
$$
The projection corresponds to the standard projection
$(x,y,z)\mapsto(x,y)$ to the coordinate $xy$-plane. It maps these
curves to the family of affine plane rational cubic curves defined by
$y^2+x^2(\tau+x)=0$ with $\tau\in \R$.

A solitary double point of the projection is not the image of any real
point of the link. It is the image of two imaginary complex conjugate
points of the complexification of the link. The preimage of the point
in the 3-space under the projection is a real line. It is disjoint from
the real part of the link, but intersects its complexification in a
couple of complex conjugate imaginary points.

In the model of the first Reidemeister move above, $(0,0)$ is the
double point of the projection for each $\tau\ne0$. If $\tau<0$, it is a
usual crossing point. Its preimage consists of two real points
$(0,0,\sqrt{-\tau})$ and $(0,0,-\sqrt{-\tau})$. If $\tau>0$, it is a
solitary double point. Its preimage consists of two imaginary conjugate
points $(0,0,i\sqrt{\tau})$ and $(0,0,-i\sqrt{\tau})$, which lie on a
real line $x=y=0$

Below, in Section \ref{s0.3}, with any solitary double point of
the projection, a local writhe equal to $\pm1$ is associated. This
is done in such a way that the local writhe of the crossing point
vanishing in the first Reidemeister move is equal to the local
writhe of the new-born solitary double point. In the case of an
algebraic knot, the sum of local writhes of all double points, both
solitary and crossings, does not depend on the choice of
projection and is invariant under rigid isotopy. This sum is the
encomplexed writhe.

\subsection{Encomplexed writhe for nonoriented and semi-oriented
links}\label{sI5} A~construction similar to the
construction of the encomplexed writhe number of an algebraic knot can
be applied to an algebraic \textit{link}. However in this case
there are two versions of the construction.

In the first of these, we define an encomplexed writhe number generalizing
the encomplexed writhe number defined above for knots. We consider a
link diagram and the sum of local writhes at solitary double
points and crossing points where the branches belong the same
connected component of the set of real points. At these crossing
point, to define a local writhe, we need orientations of the
branches. As above, we choose an orientation on each of the
components. If we make another choice, at a crossing point for
which the branches belong the same component, either both orientations
change or none. Hence the local writhe numbers at crossing points
of this kind do not depend on the choice. In Section \ref{s0.1}
below, we prove that the whole sum of local writhes over crossing
points of this kind and solitary double points does not depend on
the projection and is invariant under rigid isotopy. We call this
sum the \textit{encomplexed writhe number} of the link $A$ and denote by
$\C w(A)$.

In the second version of the construction, we consider a real
algebraic link which is equipped with an orientation of the set of
real points, use these orientations to define local writhe numbers
at all crossing points and sum the local writhe numbers over
all crossing points and all solitary double points. The result is
called the \textit{encomplexed writhe number of an oriented real algebraic
link}. This encomplexed writhe number does not change when the
orientation reverses. An orientation considered up to reversing
is called a \textit{semi-orientation}. Thus the encomplexed writhe number
depends only on the semi-orientation of the link.

The (semi-)orientation may be an artificial extra structure, but
it may also appear in a natural way, say, as a complex
orientation, if the set of real points divides the set of real
points, see \cite{R}. In fact, the complex orientation is defined
up to reversing, so it is indeed a semi-orientation. Another important
class of semi-oriented algebraic links appears as transversal
intersections of two real algebraic surfaces of degrees $p$ and
$q$ with $p\equiv q\mod2$.

The encomplexed writhe number of (semi-)oriented real algebraic link
differs from the encomplexed writhe number of the same link without
orientation by the sum of all pairwise linking numbers of the
components multiplied by $2$: let $A$ be a real algebraic link, let
$\bar A$ be the same link equipped with an orientation of its set
of real points and $\bar A_1,\dots,\bar A_n$ the (oriented)
connected components of this set, then
$$
\C w(\bar A)=\C w(A)+2\sum_{1\le i\le j\le n}
\operatorname{lk}(\bar A_i,\bar A_j).
$$

\subsection{Encomplexed writhe and framings}\label{sI6}
In the case of a knot, the encomplexed writhe
number defines a natural class of framings, since
homotopy classes of framings are enumerated by their self-linking
numbers and we can choose the framing having the self-linking
number equal to the algebraic encomplexed writhe number.
I do not know any direct construction of this framing.
Moreover, there seems to be a reason for the absence of such a
construction. In the case of links, the construction above gives a
single number, while framings are enumerated by sequences of
numbers with entries corresponding to components.

\section{Real algebraic projective links}\label{s0.1}

Let $A$ be a
nonsingular real algebraic curve in $3$-dimensional projective space.
Then the set $\R A$ of its real points is a smooth closed $1$-dimensional
submanifold of $\R P^3$, i.e., a smooth projective link. The set $\C A$
of its complex points is a smooth complex $1$-dimensional submanifold of
$\C P^3$.

Let $c$ be a point of $\R P^3$. Consider the projection
$p_c\colon\C P^3\sminus c\to \C P^2$ from $c$. Assume that $c$ is such that
the restriction to $\C A$ of $p_c$ is generic. This means
that it is an immersion without triple points and at each double point
the images of the branches have distinct tangent lines. It
follows from well-known theorems that those $c$'s for which this is the
case form an open dense subset of $\R P^3$ (in fact, it is the
complement of a $2$-dimensional subvariety).

The real part $p_c(\C A)\cap\R P^2$ of the image consists of the image
$p_c(\R A)$ of the real part and, maybe, several solitary points, which
are double points of $p_c(\C A)$.

\subsection{The local writhe of a crossing}\label{s0.2}There is a purely
topological construction which assigns a local writhe equal to $\pm1$
to a crossing belonging to the image of only one component of $\R A$.
This construction is well-known in the case of classical knots. Here is
its projective version. I borrow it from Drobotukhina's paper \cite{Dr}
on the generalization of Kauffman brackets to links in projective space.

\begin{figure}[h]
\centerline{\includegraphics[bb=0 0 125 70,clip]{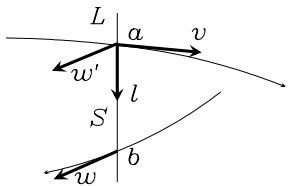}}
\caption{Construction of the frame $v$, $l$, $w'$.}
\label{f2}
\end{figure}

Let $K$ be a smooth connected one-dimensional submanifold of $\R P^3$,
and $c$ be a point of $\R P^3\sminus K$. Let $x$ be a generic double
point of the projection $p_c(K)\subset \R P^2$ and $L\subset \R P^3$ be
the line which is the preimage of $x$ under the projection. Denote by
$a$ and $b$ the points of $L\cap \R P^3$.
The points $a$ and $b$ divide the line $L$ into two segments.
Choose one of them and denote it by $S$. Choose an orientation of $K$.
Let $v$ and $w$ be tangent vectors of $K$ at $a$ and $b$ respectively
directed along the selected orientation of $K$.

Let $l$ be a vector tangent to $L$ at $a$ and directed inside $S$.
Let $w'$ be a vector at $a$ such that it is tangent to the plane
containing $L$ and $w$ and is directed to the same side of $S$ as $w$
(in an affine part of the plane containing $S$ and $w$). See
Figure~\ref{f2}.
The triple $v$, $l$, $w'$ is a base of the tangent space $T_a\R P^3$.
Define the local writhe of $x$ to be the value taken by the orientation
of $\R P^3$ on this frame.

The construction of the local writhe of $x$ contains several choices.
Here is a proof that the result does not depend on them.

We have chosen an orientation of $K$. Had the opposite orientation
been selected, then $v$ and $w'$ would be replaced by the opposite
vectors $-v$ and $-w'$. This would not change the result, since
$-v$, $l$, $-w'$ defines the same orientation as $v$, $l$, $w'$.

We have chosen the segment $S$. If the other half of $L$ was selected,
then $l$ and $w'$ would be replaced by the opposite vectors. But
$v$, $-l$, $-w'$ defines the same orientation as $v$, $l$, $w'$.

The construction depends on the order of points $a$ and $b$. The other
choice (with the same choice of the orientation of $K$ and segment $S$)
gives a triple of vectors at $b$. It can be moved continuously
without degeneration along $S$ into the triple $w'$, $-l$, $v$, which defines
the same orientation as $v$, $l$, $w'$. \qed

\subsection{Local writhe of a solitary double point}\label{s0.3}
Let $A$, $c$, and $p_c$ be as in the beginning of Section \ref{s0.1}
and let $s\in\R P^2$ be a solitary double point of $p_c$. Here is a
construction assigning $\pm1$ to $s$. I will also call the
result a \textit{local writhe} of $s$.

Denote the preimage of $s$ under $p_c$ by $L$. This is a real line in
$\R P^3$ connecting $c$ and $s$. It intersects $\C A$ in two imaginary
complex conjugate points, say, $a$ and $b$. Since $a$ and $b$ are
conjugate, they belong to different components of $\C L\sminus\R L$.

Choose one of the common points of $\C A$ and $\C L$, say, $a$. The
natural orientation of the component of $\C L\sminus\R L$ defined
by the complex structure of $\C L$ induces an orientation on $\R L$ as on
the boundary of its closure. The image under $p_c$ of the local branch
of $\C A$ passing through $a$ intersects the plane of the projection
$\R P^2$ transversally at $s$. Take the local orientation of the plane
of projection such that the local intersection number of the plane and
the image of the branch of $\C A$ is~$+1$.

Thus the choice of one of two points of $\C A\cap\C L$ defines an
orientation of $\R L$ and a local orientation of the plane of
projection $\R P^2$ (we can speak only of a local orientation of
$\R P^2$, since the whole $\R P^2$ is not orientable). The plane of
projection intersects\footnote{We may think on the plane of projection
as embedded into $\R P^3$. If you would like to think on it as on the
set of lines of $\R P^3$ passing through $c$, please, identify it in a
natural way with any real projective plane contained in $\R P^3$ and
disjoint from $c$. All such embeddings $\R P^2\to\R P^3$ are isotopic.}
transversally $\R L$ in $s$. The local orientation of the plane,
the orientation of $\R L$ and the orientation of the ambient $\R P^3$
determine the intersection number. This is the local writhe.

It does not depend on the choice of $a$. Indeed, if one chooses
$b$ instead, then both the orientation of $\R L$ and the local
orientation of $\R P^2$ would be reversed. The orientation of $\R
L$ would be reversed, because $\R L$ inherits opposite
orientations from the different halves of $\C L\sminus\R L$. The local
orientation of $\R P^2$ would be reversed, because the complex
conjugation involution $\conj\colon\C P^2\to\C P^2$ preserves the complex
orientation of $\C P^2$, preserves $\R P^2$ (point-wise) and maps one
of the branches of $p_c(\C A)$ at $s$ to the other reversing its
complex orientation.

\subsection{Encomplexed writhe and its invariance}\label{s0.4}
Now for any real algebraic projective link $A$, choose a point
$c\in\R P^3$ such that the projection of $A$ from $c$ is generic
and sum the writhes of all crossing points of the projection
belonging to the image of only one component of $\R A$ and the
writhes of all solitary double points. This sum is called the
\textit{encomplexed writhe number of $A$}.

I have to show that it does not depend on the choice of
projection. The proof given below proves more: the sum is
invariant under \textit{rigid isotopy} of $A$. By rigid isotopy we
mean an isotopy consisting of nonsingular real algebraic curves.
The effect of a movement of $c$ on the projection can be achieved
by a rigid isotopy defined by a path in the group of projective
transformations of $\R P^3$. Therefore the following theorem
implies both the independence of the encomplexed writhe number
from the choice of projection and its invariance under rigid
isotopy.

\begin{theom}\label{mainth} For any two rigidly isotopic real algebraic
projective links $A_1$ and $A_2$ whose projections from the
same point $c\in\R P^3$ are generic, the encomplexed
writhe
numbers of $A_1$ and $A_2$ defined via $c$ are equal.
\end{theom}

This theorem is proved in Section \ref{s0.6}.

\begin{cor1}\label{cor1}
The encomplexed writhe number of a real
algebraic projective link does not depend on the choice of the projection
involved in its definition.
\end{cor1}

\begin{proof}[Proof of \ref{cor1}]
A projection depends only on the
center from which it is done. The effect on the projection of a movement
of the center can be achieved by a rigid isotopy defined by a path in
the group of projective transformations of $\R P^3$.
\end{proof}

Thus the encomplexed writhe number
is a characteristic of a real algebraic link.

\begin{cor2}\label{cor2} The encomplexed writhe number of a real
algebraic projective link is invariant under rigid isotopy.\qed
\end{cor2}

\begin{figure}[htb]
\centerline{\includegraphics[bb=0 0 264 101,clip]{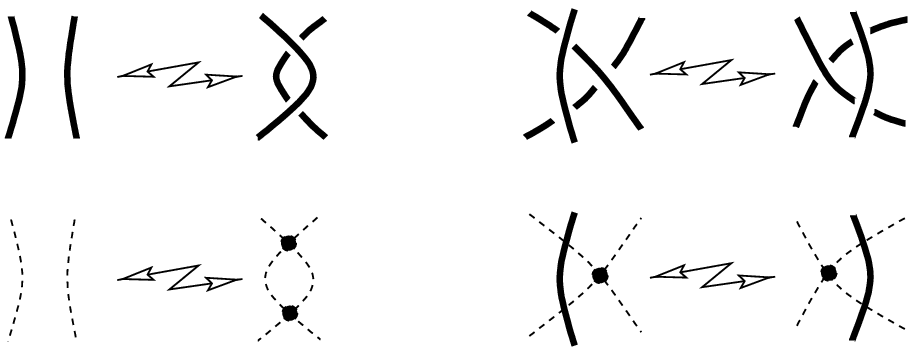}}
\caption{}
\label{o23}
\end{figure}

\subsection{Algebraic counterparts of Reidemeister moves}\label{s0.5}
As in the purely topological situation of an isotopy of a
classical link, a generic rigid isotopy of a real algebraic link may be
decomposed into a composition of rigid isotopies, each of which involves
a single local standard move of the projection. There are $5$ local
standard moves. They are similar to the Reidemeister moves. The first
of these $5$ moves is shown in the right hand half of Figure~\ref{f1}.
The other moves are shown in Figure \ref{o23}.
The first two of these coincide with the second and third Reidemeister moves.
The fourth move is similar to the second Reidemeister move: also two double
points of projection come to each other and disappear. However the
double points are solitary. The fifth move is similar to the third
Reidemeister move: a triple point also appears for a moment. But at
this triple point only one branch is real, the other two are imaginary
conjugate to each other. In this move a solitary double point
traverses a real branch.

\subsection{Reduction of Theorem \ref{mainth} to Lemmas}\label{s0.6}
To prove Theorem \ref{mainth}, first replace the rigid isotopy by
a generic one and then decompose the latter into local moves
described above, in Section \ref{s0.5}. Only in the first, fourth
and fifth moves solitary double points are involved. The
invariance under the second and the third move follows from the
well-known fact of knot theory that the topological writhe is
invariant under the second and third Reidemeister moves. Cf.~\cite{Dr}.
Thus the following three lemmas imply Theorem \ref{mainth}.

\begin{lem}\label{lem1} In the fifth move the writhe of the solitary
point does not change.
\end{lem}

\begin{lem}\label{lem2} In the fourth move the writhes of the vanishing
solitary points are opposite.
\end{lem}

\begin{lem}\label{lem3} In the first move the writhe of
vanishing crossing point is equal to the writhe of the new-born
solitary point.
\end{lem}

\subsection{Proof of Lemmas \ref{lem1} and \ref{lem2}}\label{s0.7}
Proof of Lemma \ref{lem3} is postponed to Section \ref{s0.9}.
Note that although Lemma \ref{lem3} is the most difficult to prove, it
is the least significant: here its only role is to justify the choice of
sign made in the definition of local writhe in solitary double point of
the projection. It is clear that the writhes of vanishing double
points involved in the first move are related, and if they were
opposite to each other, then the definition of
the encomplexed writhe number
should be changed, but would not be destroyed irrecoverably.

\begin{proof}[Proof of Lemma \ref{lem1}] This is obvious. Indeed, the
real branch of the projection does not interact with the imaginary
branches, it just passes through their intersection point.
\end{proof}

\begin{proof}[Proof of Lemma \ref{lem2}] At the moment of the fourth move
take a small ball $B$ in the complex projective plane centered in the
solitary self-tangency point of the projection of the curve. Its
intersection with the projection of the complex point set of the curve
consists of two smoothly embedded disks tangent to each other and to
the disk $B\cap\R P^2$. Under the move each of the disks experiences a
diffeotopy. Before and after the move the intersection the curve with
$B$ is the union of the two disks meeting each other transversally in
two points, but before the move the disks do not intersect $\R P^2$,
while after the move they intersect $\R P^2$ in their common points.

To calculate the writhe at both
vanishing solitary double points, let us select the same imaginary
branch of the projection of the curve passing through the points.
This means that we select one of the disks described above.
The sum of the local intersection numbers of this disk (equipped with
the complex orientation) and $B\cap\R P^2$ (equipped with some
orientation) is zero since under the fourth move the intersection
disappears, while in the boundary of $B$ no intersection happens.

Therefore the local orientations of the projective plane in the
vanishing solitary double points defined by this branch define opposite
orientations of $B\cap\R P^2$. (Recall that the local orientations are
distinguished by the condition that the local intersection numbers are
positive.)

On the other hand, under the move the preimages of the vanishing
solitary double points come to each other up to coincidence at the
moment of the move and their orientations defined by the choice of the
same imaginary branch are carried to the same orientation of the
preimage of the point of solitary self-tangency. Indeed, the preimages
are real lines and points of intersection of their complexifications
with the selected imaginary branch of the curve also come to the same
position. Therefore the halves of the complexifications containing the
points come to coincidence, as well as the orientations defined by the
halves on the real lines.

It follows that the intersection numbers of $B$ with the preimages of
the vanishing solitary double points equipped with these orientations
are equal. Since the local orientations of the projective plane in the
vanishing solitary double points define distinct orientations of
$B\cap\R P^2$, the writhes are opposite to each other.\end{proof}

\subsection{Proof of Lemma \ref{lem3}}\label{s0.9} It is sufficient to
consider the model family of curves described in Section \ref{SI3}.
Recall that the curves of this family are defined by the following system
of equations
$$
\left\{\begin{aligned}
xz+y&=0,\\
x+z^2+\tau&=0,
\end{aligned}\right.
$$
where $\tau$ is the parameter of the deformation. These curves admit a rational
parametrization
$$
\left\{\begin{aligned} x&=-t^2-\tau,\\
y&=-t(t^2+\tau),\\
z&=-t.
\end{aligned}\right.
$$
The projection corresponds to the standard projection
$(x,y,z)\mapsto(x,y)$ to the coordinate $xy$-plane. It maps these
curves to the family of affine plane rational cubic curves defined by
$y^2+x^2(\tau+x)=0$ with $\tau\in \R$.

We must prove that the local writhe at $(0,0)$ for $\tau<0$ coincides
with the local writhe at $(0,0)$ for $\tau>0$.

\begin{figure}[ht]
\centerline{\includegraphics[bb= 0 0 338 142,clip]{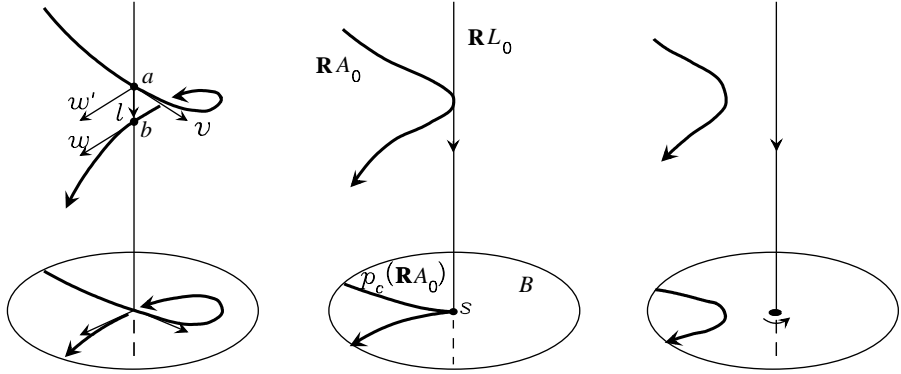}}
\caption{Real algebraic version of the first Reidemeister move.}
\label{f3}
\end{figure}

Let us calculate the local writhe for $\tau<0$. Denote $\sqrt{-\tau}$ by
$\rho$. The preimage of $(0,0)$ consists of points $a=(0,0,\rho)$ and
$b=(0,0,-\rho)$ corresponding to the values $-\rho$ and $\rho$ of $t$,
respectively, see Figure \ref{f3}. The tangent vectors to the curve at
these points are $v=(2\rho,-2\rho^2,-1)$ and $w=(-2\rho,-2\rho^2,-1)$.
The vector $l$ connecting $a$ and $b$ is $(0,0,-2\rho)$. By
definition, the writhe is the value taken by the orientation of $\R
P^3$ on the frame $v$, $l$, $w'$. This value is equal to the value of this
orientation on the frame $(1,0,0)$, $(0,1,0)$, $(0,0,1)$ multiplied by the
sign of
$$
\det\begin{pmatrix}2\rho & -2\rho^2 & -1\\
0 & 0 &
-2\rho\\
-2\rho&-2\rho^2&-1 \end{pmatrix}= -16\rho^4<0.
$$

Let us calculate the local writhe for $\tau>0$. Denote $\sqrt{\tau}$
by $\rho$. The preimage of $(0,0)$ consists of points
$a'=(0,0,i\rho)$ and $b'=(0,0,-i\rho)$ corresponding to the values
$- i\rho$ and $i\rho$ of $t$. Choose the branch which passes
through $a'$. It belongs to the upper half of the line $x=y=0$,
which induces the positive orientation of the real part directed
along $(0,0,1)$. At $a'$ the branch of the curve has tangent
vector $v=(2i\rho,2\rho^2,-1)$ and the real basis consisting of
$v$ and $iv=(-2\rho,2i\rho^2,-i)$ positively oriented with respect
to the complex orientation of this branch. The projection maps
this basis to the positively oriented basis $(2i\rho,2\rho^2)$,
$(-2\rho,2i\rho^2)$ of the projection of the branch. The
intersection number of this projection and $\R^2$ in $\C^2$ is the
sign of
$$
\det\begin{pmatrix} 0& 2\rho& 2\rho^2 &0\\
-2\rho&0 &0&2\rho^2 \\
1 & 0 & 0 & 0\\
0 & 0 &1 & 0
\end{pmatrix}=-4\rho^3<0.
$$
Hence the orientation of $\R^2$ such that its local intersection
number with the selected branch of the projection does not coincide
with the orientation defined by the standard basis. The
intersection number of the line $x=y=0$ with the standard
orientation and the $xy$-plane with the standard orientation is the
value of the orientation of the ambient space $\R^3$ taken on the
standard basis $(1,0,0)$, $(0,1,0)$, $(0,0,1)$. Therefore the
local writhe is opposite to this value. \qed

\begin{rem}\label{rem}
There is a more conceptual proof of Lemma \ref{lem3}. It is based
on a local version of the Rokhlin Complex Orientation Formula, see
\cite{R} and \cite{Viro New pr.}. In fact, the original proof was
done in that way. However, the Complex Orientation Formula is more
complicated than the calculation above.
\end{rem}

\subsection{Encomplexed writhe of an algebraic link as a Vassiliev
invariant of degree one}\label{sI4} To speak about Vassiliev
invariants, we need to fix a connected family of curves, in which
links under consideration comprise the complement to a
hypersurface. In the case of classical knots one could include all
knots in such a family by adjoining knots with self-intersections
and other singularities. A singular knot is a right equivalence
class of a smooth map of the circle to the space (recall that two
maps from a circle are right equivalent if one of them is a
composition of a self-diffeomorphism of the circle with the other
one).

In the case of real algebraic knots, such a family including all
real algebraic knots does not exist. Even the space of complex
curves in the three-dimensional projective space consists of
infinitely many components. It is impossible to change the homology
class realized by the set of complex points of an algebraic curve
in $\C P^3$ by a continuous deformation. Recall that the homology
class belongs to the group $H_2(\C P^3)=\mathbb{Z}$ and is a positive
multiple $d[\C P^1] $ of the natural generator of
$[\C P^1]\in H_2(\C P^3)$ realized by a line. The coefficient $d$ is
called the \textit{order} of the curve. The genus is another
numerical characteristic of a complex curve which takes the same
value for all nonsingular curves in any irreducible family. As is
well known, the nonsingular complex curves of given order and
genus in three-dimensional projective space are parametrized by
a finite union of quasi-projective varieties. For each of these
varieties, one can try to build a separate theory of Vassiliev
invariants on a class of nonsingular real algebraic curves whose
complexifications are parametrized by points of this variety. (A
similar phenomenon takes place in topology: links with different
numbers of components cannot be included into a reasonable
connected family, and therefore for each number of components there
is a separate theory of Vassiliev invariants.)

Among the varieties of algebraic curves in three-dimensional
projective space, there are two special families: for each natural
number $d$ there is an irreducible variety of rational curves of
order $d$ (recall that a an algebraic curve is called rational if
it admits an algebraic parametrization by a line), and for each
pair of natural numbers $p$ and $q$ there is an irreducible variety
of curves which can be presented as intersection of surfaces of
degrees $p$ and $q$.

In the class of real algebraic rational curves of order $d$,
singular curves comprise a \textit{discriminant hypersurface} in
which a generic point is a rational curve such that it has exactly
one singular point and this point is an ordinary double point. An
ordinary double point may be of one of the following two types:
either it is an intersection point of two real branches, or two
imaginary conjugate branches.

Any two real algebraic rational nonsingular curves of order $d$
can be connected by a path in the space of real rational curves of
degree $d$ that intersects the discriminant hypersurface only
transversally at a finite number of generic points. Such a path can
be regarded as a deformation of a curve to the other one. When
it intersects the discriminant hypersurface at a point, which is a
curve with singularity on real branches, the set of real points of
the curve behaves as in classical knot theory: two pieces of
the set of real points come to each other and pass through each
other. As in classical knot theory, at the moment of
intersection, the generic projection of the curve experiences an
isotopy. Nothing happens besides that one crossing point becomes
for a moment the image of a double point and then changes back into
a crossing point, but with the opposite writhe. When the path
intersects the discriminant hypersurface at a point, which is a
curve with singularity on imaginary branches, two complex
conjugate imaginary branches pass through each other. At the
moment of passing, they intersect in a real isolated double point.
At this moment the set of real points of a generic projection
experiences an isotopy. No event happens besides that a solitary
double point becomes for a moment the image of a solitary real
double point of the curve and then changes back into an ordinary
solitary double point of the projection (which is not the image of
a real point of the knot), but with the opposite writhe number.

It is clear that the encomplexed writhe number of an algebraic curve
changes under a modification of each of these kinds by $\pm2$,
with the sign depending only on the local structure of the
modification near the double point. This means that \textit{the
encomplexed writhe number on the family of real rational curves under
consideration is a Vassiliev invariant of degree $1$}.

This is true also for any space of nonsingular real algebraic
curves that can
be included into a connected family of real algebraic curves by
adjoining a hypersurface, penetration through which at a generic
point looks as in the family of rational curves described above.

There are many families of this kind besides the families of
rational knots. However, in many families of algebraic curves a
transversal penetration through the discriminant hypersurface in a
generic point looks differently. In particular, for intersections
of two surfaces it is a Morse modification of the real part of the
curve. At the moment, the old double points of the projection, both
solitary and crossing, do not change. An additional double point
appears just for a moment. However the division of
crossing points to self-crossing points of a single component and
crossing points of different components may change. Therefore the
encomplexed writhe number changes in a complicated way. If the degrees
of the surfaces defining the curve are of the same parity, the
real part of the curve has a natural semi-orientation. The Morse
modification respects this semi-orientation. Therefore the
encomplexed writhe number of the semi-oriented curve does not change.

\begin{theom}\label{vanishsl}
The encomplexed writhe number of any nonsingular semi-oriented real
algebraic link which is a transversal intersection of two real
algebraic surfaces whose degrees are of the same parity is zero.
\end{theom}

\begin{proof} Any two nonsingular real curves of the type under
consideration can be connected by a path as above. Hence their
self-linking numbers coincide. On the other hand, it
is easy to construct, for any pair of natural numbers $p$ and $q$
of the same parity, a pair of
nonsingular real algebraic surfaces of degrees $p$ and $q$
transversal to each other in
three-dimensional projective space such that their intersection
has zero self-linking number.
\end{proof}

In contrast to this vanishing result, one can prove that the
\textit{encomplexed writhe
number of a real algebraic rational knots of degree
$d$ can take any value in the interval between $-(d-1)(d-2)/2$
and $(d-1)(d-2)/2$ including these limits and congruent to
them modulo $2$}.\qed

\section{Generalizations}\label{s1}

\subsection{The case of an algebraic link with imaginary
singularities}\label{s0.10}
The same construction may be applied to
real algebraic curves in $\R P^3$ having singular imaginary
points, but no real singularities. In the construction we can
eliminate projections from the points such that some singular
point is projected from them to a real point. Indeed, for any
imaginary point there exists only one real line passing through it
(the line connecting the point with its complex conjugate), thus
we have to exclude a finite number of real lines.

This gives a generalization of encomplexed writhe numbers with the same
properties: it is invariant with respect to rigid isotopies (i.e.,
isotopies made of curves from this class), and is multiplied by
$-1$ under a mirror reflection.

\subsection{Real algebraic links in the sphere}\label{srS} The construction
of this paper can be applied to algebraic links in the sphere
$S^3$. Although from the viewpoint of knot theory this is the
most classical case, from the viewpoint of algebraic geometry the
case of curves in the projective space is simpler. The
three-dimensional sphere $S^3$ is a real algebraic variety. It is
a quadric in four-dimensional real affine space. The
stereographic projection is a birational isomorphism of $S^3$ onto
$\R P^3$. It defines a diffeomorphism between the complement of
the center of the projection in $S^3$ and a real affine space.

Given a real algebraic link in $S^3$, one may choose a real point of
$S^3$ from the complement of the link and project the link from this
point to an affine space. Then include the affine space into the
projective space and apply the construction above. The image has no
real singular points, therefore we can use the result of the previous
section.

This construction blows up the center of projection, making a real
projective plane out of it, and maps the complement to the center
of the projection in the set of real points of the sphere
isomorphically onto the complement of the projective plane. In the
imaginary domain, it contracts each generatrix of the cone which
is the intersection of the sphere with its tangent plane at the
center of projection. The image of the cone is an imaginary
quadric curve contained in the projective plane which appeared as
the result of blowing up of the central point.

\subsection{Other generalizations}\label{sGeneralizations} It is difficult
to survey all possible generalizations. Here I indicate only two
directions.

First, consider the most straightforward generalization. Let $L$ be a
nonsingular real algebraic $(2k-1)$-dimensional subvariety in the
projective space of dimension $4k-1$. Its generic projection to $\R
P^{4k-2}$ has only ordinary double points. At each double point either
both branches of the image are real or they are imaginary complex
conjugate. If the set of real points is orientable, then one can repeat
everything with obvious changes and obtain a
definition of a numeric invariant generalizing
the encomplexed writhe number defined above.

Let $M$ be a nonsingular three-dimensional real algebraic variety with
oriented set of real points equipped with a real algebraic fibration
over a real algebraic surface $F$ with fiber a projective line. There
is a construction which assigns to a real algebraic link (i.e., a
nonsingular real algebraic curve in $M$) with a generic projection to
$F$ an integer, which is invariant under rigid isotopy, is multiplied by
$-1$ under the orientation reversal in $M$ and is a Vassiliev
invariant of degree $1$. This construction is similar to the one
presented above, but uses, instead of the projection to $\R P^2$,
an algebraic version of Turaev's shadow descriptions of links \cite{T}.

\subsection{Not only writhe can be encomplexed}\label{sW}
Here we discuss only one example. However it can be easily
generalized. Consider immersions  of the sphere $S^{2n}$ to $\R^{4n}$.
Up to regular homotopy (i.e., a homotopy consisting of immersions
whose differentials also comprise a homotopy), an immersion
$S^{2n}\to\R^{4n}$ is defined by its Smale invariant \cite{Smale},
which is an element of $\pi_{2n}(V_{4n,2n})=\mathbb{Z}$. For a
generic immersion, it can be expressed as the sum of local
self-intersection numbers over all double points of the immersion,
see \cite{Smale}.

Let us encomplex the Smale invariant. For this, first, we have to
consider a real algebraic counterpart for the notion of generic
immersion $S^{2n}\to\R^{4n}$. The identification is defined via
the universal covering $\R^{4n}\to(S^1)^{4n}$. Replace Euclidean
space $\R^{4n}$ by torus $(S^1)^{4n}$, which has the advantage of
being compact. The classification of immersions
$S^{2n}\to(S^1)^{4n}$ up to regular homotopy coincides with the
Smale classification of immersions $S^{2n}\to\R^{4n}$. The sphere
$S^{2n}$ is the real part of a quadric projective hypersurface.
The torus $(S^1)^{4n}$ is the real part of a complex Abelian
variety. Consider real regular maps of the quadric to the Abelian
variety. A generic map defines an immersion both for the complex
and real parts. The only singularities are transversal double
points. Double points in the real part of the target variety are
of two kinds. At a double point of the first kind two sheets of
the image of $S^{2n}$ meet. At a double point of the second kind
the images two complex conjugate sheets of the complexification of
$S^{2n}$ meet. The Smale invariant is the sum of the local
intersection numbers over the double points of the first kind. One
can extend the definition of the local intersection number to the
double points of the second kind in such a way that the total sum
of the local intersection numbers over double points of both kinds
would be invariant under continuous deformations of regular maps.

This total sum is the \textit{encomplexed Smale invariant}. Notice
that it is, in a sense, more invariant than the original Smale
invariant. The Smale invariant may change under homotopy, it is
invariant only under regular homotopy. The encomplexed Smale
invariant does not change under a homotopy in the class of regular
maps, which corresponds to the class of all continuous
maps.

\end{document}